\documentclass[11pt]{article}
\usepackage{amssymb}
\usepackage{graphicx}
\usepackage{amsmath}
\usepackage{multirow}
\newtheorem{theorem}{Theorem}

\newtheorem{corollary}[theorem]{Corollary}

\newtheorem{lemma}[theorem]{Lemma}

\newtheorem{proposition}[theorem]{Proposition}
\newtheorem{remark}[theorem]{Remark}

\date{}
\begin{document}
\title{{\bf On the Arnold's conjecture on hyperbolic homogeneous polynomials}}
\author{
Adriana Ortiz
Rodr\'{\i}guez \thanks{
Work partially supported by DGAPA-UNAM grant PAPIIT-IN108112 and N103010} and
Federico S\'anchez-Bringas \thanks{
Work partially supported by DGAPA-UNAM grant PAPIIT-IN110803} }
\maketitle
\begin{abstract}
The Hessian Topology is a subject with interesting relations
with some classical problems of analysis and geometry
\cite{arn1}, \cite{panov}, \cite{arn2}.
In this article we prove a conjecture on this subject stated by V.I. Arnold in \cite{geoast} and  \cite{arn3}, 
concerning the number of connected components of hyperbolic homogeneous polynomials of degree $n$. The proof is
constructive and provides models. 
Our approach uses index properties at isolated singularities of
hyperbolic quadratic differential forms and combinatorial properties of recurrent functions.
\end{abstract}


\noindent {\small {\it MS classification: 53A15; 53E05; 53A05.}}

\section{Introduction}
A well known classification of the points on a smooth surface in $\mathbb R^3$ is
given in terms of the contact
of maximal order of the tangent lines with the surface at each point.
A point $p$ of a surface
is {\it elliptic} if all lines
tangent to the surface at $p$ have a contact of order
$2$ with the surface at that point. It is {\it hyperbolic} if
there are exactly two straight
lines
having a contact
of order at least $3$ with the
surface at that point. These lines are known as {\it asymptotic lines}.
A point $p$ is {\it parabolic} if it has exactly one asymptotic
line. It is possible that all the tangent lines at a point be asymptotic lines. In this case the point is named
a {\it degenerate parabolic} point.
The concept of  generic surface can be stated in terms of this type of contact in such a way that
a generic surface has the following structure:
The sets of elliptic and hyperbolic points form a union of disjoint domains on the surface whose boundary is a smooth curve
constituted by parabolic points and referred to as the {\it parabolic curve} of the surface, \cite{lan}.

Let us now consider $H^n[x,y] \subset \mathbb R[x,y]$ the set of real homogeneous
polynomials of degree $n\geq 1\,$ in two variables. The graph of any  $f \in H^n[x,y]$ contains the origin
of $\mathbb R^3$. The polynomial $f$
is called {\it hyperbolic {\rm(}elliptic{\rm)}}
if its graph is a surface with only hyperbolic (elliptic) points off the origin.
The subset of $H^n[x,y]$ constituted by hyperbolic polynomials is a topological subspace of $\mathbb R[x,y]$
denoted by $Hyp(n)$.
The connectedness of this space has been studied 
as part of the subject known as the Hessian Topology introduced in \cite{arn1}, \cite{panov}, \cite{arn2} and named 
by V. I. Arnold in \cite{geoast} and \cite{arn3} (problems 2000-1, 2000-2, 2001-1, 2002-1). In fact, in reference \cite{geoast} 
it is shown that
this property of the space depends on the degree of the polynomials that constitute it. That is,
$Hyp(3)$ and $Hyp(4)$ are connected subspaces meanwhile $Hyp(6)$ is a disconnected one.
According to this, V.I. Arnold stated the following conjecture \cite{geoast}, p.1067 and \cite{arn3}, p.139:

\begin{description}
\item[ ] {\it ``The number of connected components of the space of hyperbolic \newline homogeneous
polynomials of degree $n$ increases as $n$ increases \quad \newline
 {\rm (}at least as a linear function of $n${\rm)}."}
\end{description}
 
Moreover, the connectedness of the space of hyperbolic functions defined as follows was also studied in \cite{geoast}.
Let $n$ be a real number and $(r, \varphi)$ polar coordinates in the real plane. 
Let us define $Hyp^{\infty}(n)$ the space of smooth functions $F: S^1 
\rightarrow \mathbb R$ such that the function 
$f(r,\varphi)=r^n F(\varphi)$ referred to as a {\it homogeneous function of degree $n$} is hyperbolic, namely, its graph is constituted by
only hyperbolic points off the origin. This space has infinitely number of connected components and is closely related with the space 
of hyperbolic homogeneous polynomials.
 
In the present article we prove this conjecture, Theorem $\ref{conjetura}$. 
In fact, we present
a constructive proof. That is, we provide a good amount of examples
of hyperbolic polynomials lying on different
connected components in order to guarantee the required increasing growth of the number 
of these components in terms of the degree.
Following $\cite{geoast}$, we consider  
the field of
asymptotic lines on the graph of the polynomial. This field of lines  
has a unique singularity 
whose 
index is a convenient invariant of the connected component.
Now, let us describe the ideas we provide in order to prove the conjecture using this approach.
First, we determine inequality $(\ref{uno})$ involving a pair of polynomials $P$ and $Q$, and prove that this semi-algebraic 
condition implies 
the topological property of index preservation of the singularities of the 
fields of asymptotic lines of $P$ and $PQ$, if they are hyperbolic and $Q$ is elliptic, Corollary $\ref{coroindices}$. 
Moreover, this theorem holds not
only for polynomials but for general hyperbolic functions. 
Second, we point out that
among the families of polynomials analyzed in $\cite{geoast}$ there is no one providing an idea of the growth of the number 
of connected components of $Hyp(n)$
in terms of $n$, 
so we define a new family, using a family ${\cal P}$ analyzed in this reference. 
Family ${\cal P}$ constituted by hyperbolic homogeneous polynomials of degree $m \in \mathbb N$  
has one and only one hyperbolic polynomial of degree $m$, for each $m> 2$. Its field of asymptotic lines has a unique 
singularity at the origin with 
index $\frac{2-m}{2}$. 
On the other hand, the above new family of polynomials satisfies the following property: for a fixed degree $n$ it contains polynomials of this degree 
isotopic to those of degree $m$ lying in family ${\cal P}$ such that $m \leq n$  and $m \equiv n\ (mod\ 2)$. 
Since the index is an invariant of the connected component they belong to different components.
In order to satisfy this property we
define its elements as product polynomials $PQ$, where $P \in {\cal P}$, $Q$ is elliptic, and satisfy together inequality $(\ref{uno})$.
Let us observe that proving this inequality is equivalent to 
solving certain combinatorial equations 
involving technical formulae and recurrent functions that are not in the literature, Theorem $\ref{condPQpar}$. 
We end the article by providing a qualitative description of the foliation of the asymptotic lines of the graphs of 
these polynomials, Corollary $\ref{foliations}$.

\section{Preliminaries}

If the surface is the graph of a real valued smooth function $f$ on the plane,
the image of the parabolic curve in the $xy$-plane under the
projection $\pi: \mathbb R^3 \rightarrow \mathbb R^2,\ (x,y,z) \mapsto (x,y)$
will be referred to as the {\it Hessian
curve of $f$}. 

The directions determined by the projections of
the asymptotic lines on the $xy$-plane under this projection 
are the solutions of the following quadratic differential equation:
$$f_{xx}(x,y)dx^2+ 2f_{xy}(x,y)dxdy
+f_{yy}(x,y)dy^2=0,$$
\noindent where the quadratic differential form on the left will be denoted by $II_f(x,y)$
and referred to as the {\it second fundamental form of $f$}. Its discriminant defined as 
$$\Delta_{II_f}= f_{xy}^2- f_{xx}f_{yy},$$ 
\noindent allows us to characterize the type of point in the graph of $f$. 
That is, $(p, f(p))$ is hyperbolic (elliptic) if  $\Delta_{II_f}(p)$ is positive (negative).
The point
is parabolic if $\Delta_{II_f}(p)=0$ and $II_f(p)$ does not vanish. We say that $(p,f(p))$ is a 
degenerate parabolic point if $II_f(p)= 0$. In this case 
we say that this quadratic form has a singularity at $p$.

Let $M$ be an orientable smooth surface and
$X$ be a differentiable field of lines tangent to $M$ with an isolated singularity at the point $p$. Take a simple
closed curve $\Gamma: [0,1] \rightarrow M$,
such that $p$ is the only singularity of $X$ in the closure of the region 
determined by $\Gamma$ containing $p$.
Thus, consider the restriction of this field of lines
to the image of $\Gamma$. Moreover,
take along the image of $\Gamma$ any differentiable vector field $Y$ without singularities, 
for instance a standard coordinate vector field.
The total change of the angle between the oriented field of lines $X$
and $Y$ after going around once $\Gamma$ in the positive sense (with respect to the orientation of the surface)
is {\it the index of the field of lines} $X$ at $p$.
This number is independent of the choise of $\Gamma$ and $Y$, moreover, it has
the form $ind_p(X)=\frac{n}{2},\ n \in \mathbb Z$  \cite{hopf}.
If a tangent field of lines on $M$ is integrable the set of its integral
curves will be referred to as {\it its foliation}.
Let $X$ and $Y$ be two integrable fields of lines on $M$.  
We say that $X$ and $Y$ are {\it topologically equivalent} if there exists 
an homeomorphism $H: M \rightarrow M$ which transforms the integral curves 
of the foliation of X into the integral curves of the foliations of $Y$ \cite{nik}. 

A quadratic differential 
$$\omega(x,y)=A(x,y)dx^2 + 2B(x,y)dxdy + C(x,y)dy^2,$$  
\noindent on the punctured $xy$-plane,  $\mathbb R^2 \setminus \{(0,0) \}$ denoted by $\mathbb R^{2^*}$,
is smooth if the coefficient functions $A,B,C: \mathbb R^{2^*} \rightarrow \mathbb R$ are smooth. If 
its discriminant $\Delta_{\omega}= B^2-AC$ at a point $p$ is positive 
we will 
say that the quadratic form is {\it hyperbolic at} $p$. The quadratic form will be called {\it hyperbolic} if it is so 
at every point of its domain.  
In the sequel we will consider smooth hyperbolic quadratic differential forms whose coefficient functions extend continuously at  
the origin with values $A(0,0)=B(0,0)=C(0,0)=0$. This defines a continuous extension of $\omega$ to the plane with the origin as a unique singularity. 
The local classification of the solution curves defined by this type 
of quadratic differential forms satisfying some generic conditions 
at the singular point have been studied by several authors, see for instance 
\cite{Bruce-Tari}, \cite{Davydov} and \cite{sotom}.

The second fundamental form of a hyperbolic homogeneous polynomial, and generally that of a hyperbolic homogeneous  
function $f$ of degree $n$ are examples of this kind of smooth hyperbolic 
quadratic differential forms. Thus, $II_{f}$
defines two asymptotic lines at each point of $\mathbb R^{2^*}$. 
Moreover, it defines two continuous asymptotic fields of lines without singularities on $\mathbb R^{2^*}$ that extend
to the origin. These fields of lines are topologically equivalent.
Therefore, their indexes at the origin coincide. Consequently, this index will be called 
{\it the index of the field of asymptotic lines at the origin},
and it will be denoted by $i_0(II_f)$.




A {\it hyperbolic isotopy} between two smooth hyperbolic quadratic 
differential forms $\omega$ and $\delta$
on $\mathbb R^{2^*}$ that extend themselves to the origin with a singularity is
a smooth map
$$\Psi: \mathbb R^{2^*} \times [0,1]  \rightarrow {\cal Q},\ \ \ (x,y,t) \mapsto \Psi_t(x,y),$$
\noindent where ${\cal Q}$ is the space of real quadratic forms on the plane and the following conditions hold:
$\Psi_0(x,y)= \omega(x,y)$, $\Psi_1(x,y)= \delta(x,y)$ and $\Psi_t(x,y)$ 
is a smooth hyperbolic quadratic differential form on 
$\mathbb R^{2^*}$ which extends at the origin with a singularity. In this case we will say that {\it $\omega$
and $\delta$ are hyperbolic isotopic}.

If the second fundamental forms of two hyperbolic homogeneous polynomials of degree $n$
are hyperbolic isotopic, in fact, they are topologically equivalent.
Therefore,  
the indexes of their fields of asymptotic lines at the origin coincide.








There is a natural application  $II$, defined on the space of hyperbolic homogeneous polynomials of degree $n$, 
whose image lies on the space of 
smooth hyperbolic quadratic differential forms, associating to each polynomial its second fundamental form.  
Given two hyperbolic homogeneous polynomials of degree $n$, $f$ and $g$ lying in the same connected component $C$, 
we have a smooth curve $\gamma: [0,1] \rightarrow C$ such that $\gamma(0)=f$ and $\gamma(1)=g$. 
Then, the application $\Psi_t=II \circ \gamma(t)$ defines a
hyperbolic isotopy between $II_f$ and $II_g$.   
Therefore, let us state the following:

\begin{proposition}\label{invariante}
If two hyperbolic homogeneous polynomials of degree $n$ lie on the same connected component, the indexes of their fields of asymptotic lines
at the origin coincide.
\end{proposition}


\section{The index of the field of asymptotic lines at the origin of a hyperbolic homogeneous polynomial}
In the following analysis we consider a polynomial $f \in \mathbb R[x,y]$   
as a Hamiltonian function with
Hamiltonian vector field $\nabla f= (f_y, -f_x)$, 
on $\mathbb R^2$. The field of Hessian matrices,
$Hess f =
\left(\begin{array}{cc}
f_{xx}& f_{xy}\\
f_{xy}& f_{yy}
\end{array}
\right)$ determines at each point $p \in \mathbb R^2$ a bilinear form. 
That is, 
$$Hess f_p:\mathbb R^2  \times \mathbb R^2 \rightarrow \mathbb R $$
$$Hess f_p(X,Y)= X (Hess f(p)) Y^t,$$
\noindent where $X, Y \in \mathbb R^2$ and the index $t$, means the transpose of the vector $Y$.

Thus, for any homogeneous polynomial $P \in H^{n}[x,y]$ 
with non-null Hessian matrix
we define
the following application:
\begin{eqnarray*}
\nabla P Hess P &:& H^m[x,y] \rightarrow H^{2n+m-4}_0[x,y],\\  
&&\ \ \ \ \ \ \ Q \mapsto \nabla P{HessP} \nabla Q^t,
\end{eqnarray*}
\noindent where  $H^{2n+m-4}_0[x,y]= H^{2n+m-4} \cup \{ 0\}$.  

\noindent A straightforward computation implies that
\begin{eqnarray*}
\nabla P{HessP} \nabla Q^t= P_{xx}P_y Q_y + P_{yy} P_x Q_x -P_{xy}(P_x Q_y + P_y Q_x).
\end{eqnarray*}

\noindent The following inequality plays an important role in the proof of Theorem $\ref{teoindices}$. 
\begin{eqnarray} \label{uno}
\nabla P{HessP} \nabla Q^t(p) \leq 0,\ \ p \in \mathbb R^2.
\end{eqnarray}
In fact, if it holds at each point $p$ of the plane it guarantees 
that the isotopy applied to the second fundamental forms of the polynomials 
considered in the proof of this theorem preserves the index of the asymptotic lines at the origin.

 Let us present an interpretation of this condition.
Endow $\mathbb R^2$ with the standard orientation and interior product $<\cdot, \cdot>$. Let us define the map:
$$GradP: \mathbb R^2 \rightarrow \mathbb R^2,\ \ (x,y) \mapsto (P_x(x,y), P_y(x,y)).$$
\noindent Since the derivative of this map at $p_0$ is $D_{p_0}Grad P= {Hess P}_{p_0}$, if we assume that  
$\nabla Q^t(p_0)$ is not in the kernel of ${HessP}_{p_0}$, 
we have that the image by this map of the level curve of the polynomial $Q$ at $p_0$ 
has tangent vector ${HessP}_{p_0} \nabla Q^t(p_0)$. 
Using the parallel translation of $\mathbb R^2$, we can suppose that this curve 
intersects the level curve of the polynomial $P$ at $p_0$.
Thus, inequality $(\ref{uno})$ holds at each point of the plane, if and only if  
the oriented angle of intersection of these curves 
at each point where the polynomial function $P$ is not singular lies on the interval  $\left[ \frac{\pi}{2}, \frac{3\pi}{2} \right]$. 
\begin{theorem}\label{teoindices}
Let $P, Q$ be homogeneous polynomials such that $P$ is hyperbolic,  $Q$ is elliptic
and the
product $\, f= P Q $ is hyperbolic. Suppose also that $Q$ is positive on $\mathbb R^{2^*}$
and inequality $(\ref{uno})$ holds at each point of the plane.
Then, $II_f$ and $II_P$ are hyperbolic isotopic.
\end{theorem}

\noindent {\bf Proof.} Let us present the proof divided in lemmas.
A straightforward computation shows the following
\begin{lemma} The discriminant  of the quadratic differential form
$Q II_p+2dPdQ$ has the following expression:
\begin{eqnarray} \label{Delt}
\Delta_{Q II_p  + 2  dP dQ}  =  - Q^2 det(Hess P)+ 4\Delta_{dP dQ}
- 2Q
(\nabla P{HessP} \nabla Q^t).
\end{eqnarray}
\end{lemma}

\begin{lemma}\label{eliminacionformacuadratica}
Suppose that  $P $ is a hyperbolic homogeneous polynomial, $Q$ is a positive polynomial on
 $\mathbb R^{2^*}$ and assume that inequality $(1)$ holds.
Then,

\noindent a) The quadratic differential form $Q II_p  + 2  dP dQ$ is hyperbolic on $\mathbb R^{2^*}$.

\noindent b) The quadratic differential forms $Q II_p  + 2  dP dQ$ and $Q II_p$ are hyperbolic isotopic.
\end{lemma}
\noindent {\bf Proof.}
a) Since $P$ is a hyperbolic polynomial and 

$$\Delta_{dP dQ}= \frac{1}{4}(P_x Q_y - P_y Q_x)^2,$$ 

\noindent then inequality $(\ref{uno})$ implies  
that the right side of equation $(\ref{Delt})$  is positive on $\mathbb R^{2^*}$, that is, $Q II_p  + 2  dP dQ$ 
is hyperbolic on $\mathbb R^{2^*}$.
\newline b) Considering the isotopy $\Psi_t(x,y)= Q II_p + 2t dP dQ (x,y),
\,t \in [0,1]$,
we can see from equation $(\ref{Delt})$ that
the discriminant of the  quadratic differential form 
$Q II_p  + 2t  dP dQ (x,y)$ is
$$\Delta_{\Psi_t}= -Q^2 det(Hess P) + t^2 
\Delta_{dPdQ} - 2t Q(\nabla P{HessP} \nabla Q^t),$$ 
which is positive on $\mathbb R^{2^*}$.
This implies that $\Psi_t(x,y)$ is a hyperbolic isotopy.
$\hfill\Box$

\medskip

In order to prove the following lemma we point out the next easy
\begin{remark} \label{isohyper}
Let $\, a, b, c $ be real numbers such that $\, a + b + c > 0 ,\,\, a > 0$ and $c \leq 0$. Then
$\, a  + tb + t^2 c > 0 \,$ for  $ \ t \in (0,1)$.
\end{remark}

\begin{lemma}\label{isoformas}
Let $\,\omega, \delta $ be two smooth quadratic differential forms on $\mathbb R^2$ vanishing at the origin.  Suppose that
$\omega$, $\omega+ \delta $ are hyperbolic and $\delta$ is non-hyperbolic at
each point of $\mathbb R^{2^*}$.
Then, the quadratic forms $\omega$ and $\omega + \delta$ are hyperbolic isotopic.
\end{lemma}
\noindent {\bf Proof.}
Let  $\,\omega(x,y)=\omega_1 dx^2 + 2 \omega_2 dx dy + \omega_3 dy^2$
and $\delta(x,y)= \delta_1 dx^2 + 2 \delta_2 dx dy + \delta_3 dy^2$ be two quadratic differential forms.
Consider the isotopy
$$ \Psi_t(x,y)= \omega(x,y) + t \delta(x,y). $$
The discriminant of
$ \Psi_t(x,y)$ is
$$   \Delta \Psi_t = \omega_2^2 -  \omega_1 \omega_3 + t (2 \omega_2\delta_2- \omega_1\delta_3 \omega_3\delta_1) 
+ t^2 (\delta_2^2 - \delta_1 \delta_3).$$

\medskip
Because $\omega$  and  $\omega +\delta$ are hyperbolic on $\mathbb R^{2^*}$,
we have that 
$$\,\omega_2^2 -  \omega_1 \omega_3 > 0\ \ {\rm and}\ \
\omega_2^2 -  \omega_1 \omega_3 +  (2 \omega_2\delta_2- \omega_1\delta_3 \omega_3\delta_1) + (\delta_2^2 -  \delta_1 \delta_3)>0,$$ 
\noindent on $\mathbb R^{2^*}$. Moreover, since $\delta_2^2 -  \delta_1 \delta_3 \leq 0$
the fact that $\Delta \Psi_t(x,y)$ is negative at each point on $\mathbb R^{2^*}$ follows from Remark \ref{isohyper}. $\hfill\Box$

\medskip
\noindent {\bf End of the proof of Theorem $\ref{teoindices}$.}
Let us take  $\, \omega = Q II_p  + 2 dPdQ$
and $\, \delta = P II_Q$. Observe that Lemma $\ref{eliminacionformacuadratica}\ a)$ ensures that $\omega$ is hyperbolic on $\mathbb R^{2^*}$.
Moreover, since $f=PQ$ is hyperbolic we  conclude, by Lemma $\ref{isoformas}$
that there exists a hyperbolic isotopy between
the quadratic forms
$II_{f}= P II_Q + Q II_P + 2 dP dQ$
and $\,  Q II_p  + 2 dP dQ$. Then, Lemma $\ref{eliminacionformacuadratica}\ b)$ implies that $II_f$ and $II_P$ are hyperbolic isotopic.
$\hfill\Box$

\begin{corollary}\label{coroindices}
Assume that $f$ and $P$ satisfy the hypothesis of Theorem $\ref{teoindices}$. Then the fields of asymptotic lines  of these polynomials 
extended to
the origin with a singularity
are topologically equivalent on $\mathbb R^2$ and their indexes at the origin coincide. 
\end{corollary}

\begin{remark} We observe that Corollary $\ref{coroindices}$ is true 
in a more general setting. Namely, 
let us consider instead of homogeneous polynomials, a pair of real valued differentiable functions $P$ and $Q$ vanishing at the origin. Suppose that    
they satisfy the hyperbolicity and ellipticity hypothesis of the statement in $\mathbb R^{2^*}$, respectively. Assume that $Q$ is positive in $\mathbb R^{2^*}$ and
the origin is a degenerate parabolic point of $f=PQ$ and $P$. Suppose that inequality $(1)$, 
stated in this case on the class of differentiable functions holds 
at each point of the plane. Then, the fields of asymptotic lines of these functions are
topologically equivalent on $\mathbb R^2$ and their indexes at the origin coincide. 
\end{remark}

\section{Proof of the conjecture}

For clearness sake we present the {\it strategy proof} of the  conjecture:
In order to determine the desired number of connected components of the hyperbolic homogenous polynomials in terms of the degree,
the goal is to find for each $n \in \mathbb N$ a big enough number of polynomials in $Hyp(n)$ whose
fields of asymptotic lines at the origin have different indexes.
For this, we first consider the well known family $\cal{P}$  
of hyperbolic
homogeneous polynomials of degree $m$, whose elements $P^{m}$ described below,
define fields of asymptotic lines at the origin of indexes $\frac{2-m}{2}, m>2$, respectively, \cite{geoast} p.1037. Observe that
each element of the family determines only one connected component of $Hyp(m)$.
On the other hand, for each $n \geq 4$ 
we consider an elliptic homogeneous polynomial $Q^{2k}$ of
degree $2k$ 
in such a way that the family of polynomials $\{P^m Q^{2k} \}$, where $1 \leq k \leq \frac{n}{2} -1$ and $n=m+2k$
is contained in $Hyp(n)$, Proposition $10$. The corresponding family of fields of asymptotic lines
have indexes $\frac{2-m}{2}$ at the origin, respectively.
That is, the contribution of $Q^{2k}$ to the index at the origin of the field of asymptotic lines of $P^m Q^{2k}$ is null.
To prove that this index property holds for this type polynomials  
we apply Corollary $\ref{coroindices}$. Namely, we need to prove that any pair of polynomials $P \in \{P^m \}$ and $Q \in \{ Q^{2k}\}$ 
satisfies condition $(\ref{uno})$,
Theorem $\ref{condPQpar}$.

Let us describe the two families of homogeneous polynomials that we will use in the proof.

\medskip
In \cite{geoast}, V. I. Arnold proves that the polynomials of degree $m$
$$ P(x,y) = r^{m-k} \mbox{Re} (x+i y)^k ,$$
where $r= \sqrt{x^2+y^2} ,\ k^2> m, \, k\leq m$ and $\,m-k\,$ is even, are hyperbolic homogeneous
polynomials
and $\,i_0 (II_P) = \frac{2-k}{2}$.
Taking, in particular $ \, k =m \geq 2$ these polynomials get the form
$$ P^m(x,y)= \sum_{j=0}^\frac{m}{2} (-1)^j \left(
\begin{array}{c}
m\\ 2j\\
\end{array} \right) x^{m-2j} y^{2j} , \quad  \mbox{ if } m \mbox{ is even},$$
and
$$ P^m(x,y)= \sum_{j=0}^\frac{m-1}{2} (-1)^j \left(
\begin{array}{c}
m\\ 2j\\
\end{array} \right) x^{m-2j} y^{2j} , \quad  \mbox{ if } m \mbox{ is odd},$$
with (in both cases)  $\,
\,i_0 (II_{P_m}) = \frac{2-m}{2}.$
\bigskip

Now, let us consider the family of polynomials
$$Q^{2k}(x,y)= (x^2 + y^2)^k,$$
\noindent where $k$ is a positive integer.

A direct computation of the discriminant of the 
form $\,II_{Q^{2k}}$ implies
\begin{proposition} \label{Qelliptic}
For $\, k\geq 1$ the polynomial $\, Q^{2k}$ is elliptic.
\end{proposition}

\begin{proposition} \label{polhyper}
Let $k, m \in \mathbb N$ such that $m>max \{2,k\}$.
Then,
the homogeneous polynomial
$ f(x,y)= P^{m}(x,y) Q^{2k}(x,y) \,$   is hyperbolic.
\end{proposition}

\noindent {\bf Proof.}
Consider a homogeneous polynomial $f(x,y)$ of degree $n$, such that in polar 
coordinates $(r, \varphi)$, $f(r,\varphi)=r^nF(\varphi)$, where $F(\varphi)$ is a trigonometric function.
The following hyperbolicity condition stated by 
V. I. Arnold in \cite{geoast} (Theorem 1 p.1031) 
guarantees that
$f$ is hyperbolic if and only if the function $F$ satisfies
\begin{equation}  \label{cond-hyper}
n^2 F^2 + n F F'' - (n-1) (F')^2 < 0.
\end{equation}  
In our case, $f(r, \varphi)= r^{m+2k} \cos(m\varphi).$
Thus, the left-hand side of this inequality 
has the expression
\begin{equation} \label{condhyper}
n^2 F^2 + n F F'' - (n-1) (F')^2 = \cos^2(m\varphi) [4k (m+k)] - m^2(m+2k-1).
\end{equation}
Note that the right-side of (\ref{condhyper}) is negative since
$$\qquad\qquad\qquad\qquad\quad 4k (m+k) < m^2(m+2k-1). \hskip 3.7cm \Box $$

\begin{theorem} \label{condPQpar}
Let  $k, m \in \mathbb Z$ such that $k \geq 0$ and  $m \geq 2$.
Then, the polynomials  $P^m, Q^{2k}$
satisfy inequality $(\ref{uno})$.
\end{theorem}
The most important part of the proof of this theorem is a consequence of some combinatorial relations 
which are not in the literature. Now, we prove them. 
We present only the case when $m$ is even because the odd case is analogous.
Let us begin by proving the formulae below.

\begin{lemma}\label{lemaABC}
Let $m\geq 2\,$ be a natural number. For each integer number $\,0\leq j\leq \frac{m-2}{2}\,$
consider the combinatorial functions 
\begin{eqnarray*}
A(j) &=&  (-1)^{j} \left[
\left(
\begin{array}{c}
m-1\\ 2j\\
\end{array}
\right) +
 \sum_{k=0}^{j-1} \left[
\left(
\begin{array}{c}
m-1\\ 2k\\
\end{array}
\right)
\left(
\begin{array}{c}
m-1\\ 2j-2k\\
\end{array}
\right) - \right.\right.\qquad\qquad\nonumber\\
&& \left.\left.
\left(
\begin{array}{c}
m-1\\ 2k+1\\
\end{array}
\right)
\left(
\begin{array}{c}
m-1\\ 2j-2k-1\\
\end{array}
\right)
\right]\right],
\end{eqnarray*}

{
\begin{eqnarray*}
B = (-1)^{\frac{m}{2}}\left(
1-m + \sum_{k=0}^{\frac{m}{2}-2}
\left(
\begin{array}{c}
m-1\\ 2k+1\\
\end{array}
\right)
\left[
\left(
\begin{array}{c}
m-1\\ 2k+2\\
\end{array}
\right)-
\left(
\begin{array}{c}
m-1\\ 2k\\
\end{array}
\right)
\right]
\right),\quad
\end{eqnarray*}}
{
\begin{eqnarray*}
C(j) &=& (-1)^{\frac{m}{2}+j-1}  \left[
- \left(
\begin{array}{c}
m-1\\ 2j-1\\
\end{array}
\right) +
 \sum_{k=0}^{\frac{m}{2}-j-1} \left[
\left(
\begin{array}{c}
m-1\\ 2k+2j\\
\end{array}
\right)
\left(
\begin{array}{c}
m-1\\ 2k+1\\
\end{array}
\right)\right.\right.\quad\nonumber\\
&& \left.\left.
 -  \left(
\begin{array}{c}
m-1\\ 2k\\
\end{array}
\right)
\left(
\begin{array}{c}
m-1\\ 2k+2j-1\\
\end{array}
\right)
\right]\right].
\end{eqnarray*}}
\end{lemma}

\noindent {\it Then, they can be reduced to the following expressions.}
\begin{eqnarray}\label{coefA}
&&A(j) = \left(
\begin{array}{c}
m-1\\ j\\
\end{array}
\right), \\ \label{coefB}
&&B = \left(
\begin{array}{c}
m-1\\ \frac{m}{2}\\
\end{array}
\right),\\ \label{coefC}
&&C(j) = \left(
\begin{array}{c}
m-1\\ j+\frac{m}{2}-1\\
\end{array}
\right).
\end{eqnarray}

\noindent {\bf Proof.}

Let us begin by proving $(\ref{coefA})$.
Using several times the formula 
\begin{eqnarray}\label{fabsorcion}
(m-k) \left(\begin{array}{c}
m\\ k\\
\end{array}\right)= m \left(\begin{array}{c}
m-1\\ k\\
\end{array}\right),
\end{eqnarray}
(derived from the absorption identity \cite{graham}), the expression (\ref{coefA})
becomes
{\small
\begin{eqnarray}\label{Asimplificada}
(-1)^{j} \left[
\left(
\begin{array}{c}
m-1\\2j\\
\end{array}
\medskip
\right) +
\sum_{k=0}^{j-1}\frac{4k-2j+1}{m}
\left(
\begin{array}{c}
m\\2k+1\\
\end{array}
\right)
\left(
\begin{array}{c}
\medskip
m\\2j-2k\\
\end{array}
\right)\right]
=
\left(
\begin{array}{c}
m-1\\j\\
\end{array}
\right)
\end{eqnarray}}
By the formula for the alternating sum of consecutive binomial coefficients,
\begin{eqnarray}\label{alternante}
(-1)^r \left(
\begin{array}{c}
\medskip
m-1\\ r\\
\end{array}
\right) = \sum_{k=0}^{r} (-1)^k
\left(
\begin{array}{c}
m\\ k\\
\end{array}
\right),
\end{eqnarray}
the expression (\ref{Asimplificada}) results
{\small
\begin{eqnarray*}
(-1)^{j} \sum_{k=0}^{j-1}\frac{4k-2j+1}{m}
\left(
\begin{array}{c}
m\\2k+1\\
\end{array}
\right)
\left(
\begin{array}{c}
\medskip
m\\2j-2k\\
\end{array}
\right)
=  \sum_{k=1}^{j} (-1)^{k+1}
\left(
\begin{array}{c}
m\\ j+k\\
\end{array}
\right).
\end{eqnarray*}}

For the following reduction it is useful to write down the sum of both 
sides of the last expression in two parts. 
For the left side sum, chose the first part as the sum from the lowest 
value of $k$ up to $\frac{j}{2}$ ($\frac{j+1}{2}$) if  $j$ is even (odd), 
and the second part containing the remaining terms.  For the right side sum,
take the first part as the sum of even terms and the second part as the sum
of odd terms. Then, by associating corresponding terms of both sides
we obtain the equation

\begin{eqnarray}\label{criel}
\sum_{k=1}^{j} (-1)^{k}
\left(
\begin{array}{c}
m\\j+k\\
\end{array}
\right)
\left[ 1 + \frac{1-2k}{m}
\left(
\begin{array}{c}
m\\ j-k+1\\
\end{array}
\right)\right] = 0.
\end{eqnarray}
We present now a proof of (\ref{criel}) based on some recurrence
relations. It was given by C. Merino L\'opez \cite{Merino}. 
\medskip 

By the alternating sum (\ref{alternante}) we have
{\scriptsize
\begin{eqnarray*}
m \sum_{k=0}^{j} (-1)^k
\left(
\begin{array}{c}
m\\j+k\\
\end{array}
\right)
= 
\left\{
\begin{array}{c}
(m-2j)\left(
\begin{array}{c}
m\\2j\\
\end{array}
\right)
- (m-j) \left(
\begin{array}{c}
m\\j\\
\end{array}
\right) \quad\mbox{{\normalsize if $j$ is even}}\\
-(m-2j)\left(
\begin{array}{c}
m\\2j\\
\end{array}
\right)
- (m-j) \left(
\begin{array}{c}
m\\j\\
\end{array}
\right) \quad\mbox{{\normalsize if $j$ is odd}}\\
\end{array}
\right.
\end{eqnarray*}}
By formula (\ref{fabsorcion}) the last equality becomes

{\scriptsize
\begin{eqnarray*}
m \sum_{k=0}^{j} (-1)^k
\left(
\begin{array}{c}
m\\j+k\\
\end{array}
\right)
= 
\left\{
\begin{array}{c}
(2j+1)\left(
\begin{array}{c}
m\\2j+1\\
\end{array}
\right)
- (j+1) \left(
\begin{array}{c}
m\\j+1\\
\end{array}
\right) \quad\mbox{{\normalsize if $j$ is even}}\\
-(m-2j)\left(
\begin{array}{c}
m\\2j\\
\end{array}
\right)
- (m-j) \left(
\begin{array}{c}
m\\j\\
\end{array}
\right) \quad\mbox{{\normalsize if $j$ is odd}}\\
\end{array}
\right.
\end{eqnarray*}}
Replacing the last equality in (\ref{criel}) we obtain
\begin{eqnarray*}
\sum_{k=1}^{j+1}(-1)^{k+1} (2k-1)
\left(
\begin{array}{c}
m\\ k+j\\
\end{array}
\right)
\left(
\begin{array}{c}
m\\ j-k+1\\
\end{array}
\right)
=  (j+1) \left(
\begin{array}{c}
m\\ j+1\\
\end{array}
\right).
\end{eqnarray*}
Denote by $\,T(m,j)\,$ the left side of the last expression. Using
the {\it Stifel identity}
{\scriptsize $\,\left(
\begin{array}{c}
m\\j\\
\end{array}
\right) = \left(
\begin{array}{c}
m-1\\j\\
\end{array}
\right) + \left(
\begin{array}{c}
m-1\\j-1\\
\end{array}
\right)$},
we verify that $\,T(m,j)\,$ satisfies
the recurrence relation
\begin{eqnarray}\label{recurrenciaT}
T(m,j) &=& T(m-1,j) + T(m-1,j-1) + \left(
\begin{array}{c}
m-1\\ j\\
\end{array}
\right)^2 + \nonumber\\
&&+ \,\, 2
\sum_{k=1}^{j}(-1)^{k}
\left(
\begin{array}{c}
m-1\\ j-k\\
\end{array}
\right)
\left(
\begin{array}{c}
m-1\\ j+k\\
\end{array}
\right).
\end{eqnarray}
Now, using again the Stifel formula for the function
\begin{eqnarray*}
F(m-1,j) = \left(
\begin{array}{c}
m-1\\ j\\
\end{array}
\right)^2 + 2
\sum_{k=1}^{j}(-1)^{k}
\left(
\begin{array}{c}
m-1\\ j-k\\
\end{array}
\right)
\left(
\begin{array}{c}
m-1\\ j+k\\
\end{array}
\right),
\end{eqnarray*}
we have that it verifies the recurrence relation
\begin{eqnarray}\label{recurrenciaF}
F(m,j) = F(m-1,j) + F(m-1,j-1).
\end{eqnarray}
Considering the Stifel identity, we remark that the expression {\scriptsize $\,\left(
\begin{array}{c}
m\\j\\
\end{array}
\right)$}\, satisfies also the recurrence relation (\ref{recurrenciaF}).
Because its initial values are the same that those of (\ref{recurrenciaF})
we conclude that  $\, F(m,j) =$ {\scriptsize$\left(
\begin{array}{c}
m\\j\\
\end{array}
\right)$}.
So, the recurrrence relation (\ref{recurrenciaT}) becomes
$$ T(m,j) - T(m-1,j) - T(m-1,j-1) =
\left(
\begin{array}{c}
m-1\\j\\
\end{array}
\right).$$
But this relation is also satisfied by
$\, (j+1) ${\scriptsize $\left(
\begin{array}{c}
m\\j+1\\
\end{array}
\right)$} and moreover, their initial values are the same.
Then $ T(m,j) = (j+1) ${\scriptsize
$\left(\begin{array}{c}
m\\j+1\\
\end{array}\right)$}. This proves (\ref{criel}).

\normalsize

\bigskip

\noindent Now, let us prove equation  $(\ref{coefC})$.
By formula (\ref{fabsorcion}) the expression $C(j)$ results
\begin{eqnarray*}
C(j) &=&
(-1)^{\frac{m}{2}+j-1} \left[
\sum_{k=0}^{\frac{m}{2}-j}\frac{m-4k-2j-1}{m}
\left(
\begin{array}{c}
m\\ 2k+1\\
\end{array}
\right)
\left(
\begin{array}{c}
m\\ 2k+2j\\
\end{array}
\right)\right].
\end{eqnarray*}
So, the expression (\ref{coefC}) becomes
{\small
\begin{eqnarray}\label{Csimplificado}
(-1)^{\frac{m}{2}+j-1} \left[
\sum_{k=0}^{\frac{m}{2}-j}\frac{m-4k-2j-1}{m}
\left(
\begin{array}{c}
m\\ 2k+1\\
\end{array}
\right)
\left(
\begin{array}{c}
m\\ 2k+2j\\
\end{array}
\right)\right]
=
\left(
\begin{array}{c}
m-1\\ j+\frac{m}{2}-1\\
\end{array}
\right)
\end{eqnarray}}

\medskip
\noindent Because $m$ is even we replace $\,\, m = 2r\,$ in both sides 
of the last expression. Moreover, we consider the change $\,\, r-j = n\,$ 
to obtain
{\small
\begin{eqnarray*}
(-1)^{2r-n-1} \left[
\sum_{k=0}^{n}\frac{(2n-4k-1)}{2r}
\left(
\begin{array}{c}
2r\\ 2k+1\\
\end{array}
\right)
\left(
\begin{array}{c}
2r\\ 2k+2r-2n\\
\end{array}
\right)\right]
=
\left(
\begin{array}{c}
2r-1\\ 2r-n-1\\
\end{array}
\right).
\end{eqnarray*}}

\medskip
\noindent
Using in the last expression the symmetry identity
\,{\scriptsize $ \left(
\begin{array}{c}
a\\ b\\
\end{array}
\right)
= \left(
\begin{array}{c}
a\\ a-b\\
\end{array}
\right)\,$} 
\cite{graham}, and 
replacing to $\, 2r\,$ by $m$ it results

\begin{eqnarray}\label{f3}
(-1)^{n} \left[
\sum_{k=0}^{n}\frac{(4k-2n+1)}{m}
\left(
\begin{array}{c}
m\\ 2k+1\\
\end{array}
\right)
\left(
\begin{array}{c}
m\\ 2n-2k\\
\end{array}
\right)\right]
=
\left(
\begin{array}{c}
m-1\\ n\\
\end{array}
\right)
\end{eqnarray}

\medskip
\noindent Note that the corresponding term $\,k=n\,$ on the left is
{\scriptsize $\left(
\begin{array}{c}
m-1\\ 2n\\
\end{array}
\right)$}. 
Since expression (\ref{f3}) becomes (\ref{Asimplificada}), equation (\ref{coefC})
is proved.
\bigskip

We conclude by proving equation $(\ref{coefB})$.
Equation $(\ref{fabsorcion})$ implies that
\begin{eqnarray*}
\left( \begin{array}{c}
m-1\\ 2k+2\\
\end{array}
\right)-
\left(
\begin{array}{c}
m-1\\ 2k\\
\end{array}
\right) =
\frac{m-4k-3}{m}
\left(
\begin{array}{c}
m\\ 2k+1\\
\end{array}
\right)
\left(
\begin{array}{c}
m\\ 2k+2\\
\end{array}
\right).
\end{eqnarray*}
Then $B(j)$ results
\begin{eqnarray*}
(-1)^{\frac{m}{2}}\left(
\sum_{k=0}^{\frac{m}{2}-1}
\frac{(m-4k-3)}{m}
\left(
\begin{array}{c}
m\\ 2k+1\\
\end{array}
\right)
\left(
\begin{array}{c}
m\\ 2k+2\\
\end{array}
\right)
\right).
\end{eqnarray*}

\noindent So, we must prove
\begin{eqnarray}\label{f2}
(-1)^{\frac{m}{2}}\left(
\sum_{k=0}^{\frac{m}{2}-1}
\frac{(m-4k-3)}{m}
\left(
\begin{array}{c}
m\\ 2k+1\\
\end{array}
\right)
\left(
\begin{array}{c}
m\\ 2k+2\\
\end{array}
\right)
\right)=
\left(
\begin{array}{c}
m-1\\ \frac{m}{2}\\
\end{array}
\right)
\end{eqnarray}
But, when we replace $\, j=1$ in (\ref{Csimplificado}),
we retrieve (\ref{f2}).
\hfill$\Box$
\bigskip

\noindent {\bf Proof of Theorem \ref{condPQpar}.}
In order to prove that inequality $(\ref{uno})$ holds for the 
polynomials $P^m$ and $Q^{2k}$ we consider the polynomial expression of
\linebreak 
$\nabla P^m {HessP} (\nabla Q^{2k})^t$ and prove the following
$$ P_x(Q_y P_{xy}  - Q_x P_{yy}) + P_y (Q_x P_{xy} - Q_y P_{xx}) = 2\,k \, m^2 (m-1) (x^2+y^2)^{k+m-2}.$$

\noindent Since $m$ is even, a straightforward computation shows that
$$Q_y P_{xy} - Q_x P_{yy} =2k (x^2+y^2)^{k-1}
\left[
\sum_{j=0}^{\frac{m}{2}-1} (-1)^j
\frac{ (m-1) m!}{(2j)! (m-2j-1)!}
x^{m-2j-1} y^{2j} \right].$$

\noindent Now, we multiply both sides of the last expression by $P_x$. The
product $\,P_x \,(Q_y P_{xy} - Q_x P_{yy})\,$ results
\medskip
\begin{eqnarray*}
 2k (x^2+y^2)^{k-1} m^2 (m-1) \left[
\sum_{j=0}^{\frac{m}{2}-1} (-1)^j
\left(
\begin{array}{c}
m-1\\ 2j\\
\end{array}
\right)
x^{m-2j-1} y^{2j} \right]^2.
\end{eqnarray*}

\noindent Developing the squared factor of the last expression we have
$$P_x \,(Q_y P_{xy} - Q_x P_{yy}) \,=\, 2\,k \,(x^2+y^2)^{k-1} m^2 (m-1) \hskip 4.3cm$$
$$ \quad\left[
x^{2m-2} \,+ \sum_{j=1}^{\frac{m}{2}-1} \left(\sum_{k=0}^j (-1)^j
\left(
\begin{array}{c}
m-1\\ 2k\\
\end{array}
\right)
\left(
\begin{array}{c}
m-1\\ 2j-2k\\
\end{array}
\right)\right)
x^{2m-2j-2} y^{2j} \,\, +  \right.$$
{\small
\begin{eqnarray}\label{ecua1}
\left.\qquad\sum_{j=1}^{\frac{m}{2}-1} 
\left(\sum_{k=0}^{\frac{m}{2}-j-1} (-1)^{\frac{m}{2}+j-1}
\left(
\begin{array}{c}
m-1\\ 2k+2j\\
\end{array}
\right)
\left(
\begin{array}{c}
m-1\\ m-2k-2\\
\end{array}
\right)\right)
x^{m-2j} y^{m+2j-2}\right]
\end{eqnarray}}

\medskip

\noindent Now, we shall compute the expression $ P_y (Q_x P_{xy} - Q_y P_{xx})$. 
After doing some elemental simplifications we have 
{\small $$Q_x P_{xy} - Q_y P_{xx} = 2k (x^2+y^2)^{k-1} 
\left[\sum_{r=0}^{\frac{m}{2}-1} (-1)^{r+1} m (m-1)
{\mbox
{\scriptsize$
\left(
\begin{array}{c}
m-1\\ 2r+1\\
\end{array}
\right)$}}
x^{m-2r-2} y^{2r+1}\right].$$}

\noindent Now, consider the product of the last expression by $P_y$. So,
the product $\, P_y \,(Q_x P_{xy} - Q_y P_{xx})\,$ results
$$2k (x^2+y^2)^{k-1} (m-1) m^2 y^2
\left[
\sum_{j=0}^{\frac{m}{2}-1} (-1)^{j+1}
\left(
\begin{array}{c}
m-1\\ 2j+1\\
\end{array}
\right)
x^{m-2j-2} y^{2j} \right]^2.$$

\noindent Developing the squared term we obtain

$$P_y \,(Q_x P_{xy} - Q_y P_{xx}) \,=\, 2\,k \,(x^2+y^2)^{k-1} m^2 (m-1)\hskip 4.2cm$$
{\small
\begin{eqnarray*}
\qquad\left[ y^{2m-2} + \sum_{j=1}^{\frac{m}{2}} \left(\sum_{k=0}^{j-1} (-1)^{j+1}
\left(
\begin{array}{c}
m-1\\ 2k+1\\
\end{array}
\right)
\left(
\begin{array}{c}
m-1\\ 2j-2k-1\\
\end{array}
\right)\right)
x^{2m-2j-2} y^{2j} \,\,+ \right.
\end{eqnarray*}}
{\small
\begin{eqnarray}\label{ecua2}
\left.\qquad\sum_{j=2}^{\frac{m}{2}-1} \left(\sum_{k=0}^{\frac{m}{2}-j} (-1)^{\frac{m}{2}+j}
\left(
\begin{array}{c}
m-1\\ 2k+2j-1\\
\end{array}
\right)
\left(
\begin{array}{c}
m-1\\ m-2k-1\\
\end{array}
\right)\right)
x^{m-2j} y^{m+2j-2}\right]
\end{eqnarray}}

 Adding the expressions  (\ref{ecua1}) and (\ref{ecua2}) we obtain
 
\begin{eqnarray}\label{suma-simpli} \nonumber
P_x \,(Q_y P_{xy} - Q_x P_{yy}) + P_y \,(Q_x P_{xy} - Q_y P_{xx}) = \hskip 4.5cm\\ \nonumber
= 2\,k \,(x^2+y^2)^{k-1} m^2 (m-1)\Bigg[ x^{2m-2}\quad\hskip 2.4cm\\ \nonumber
 + \bigg(\sum_{j=1}^{\frac{m}{2}-1} A(j) x^{2m-2j-2} y^{2j}\bigg) + B x^{m-2} y^{m}
 \hskip 2.4cm \\ 
+ \left. \bigg(\sum_{j=2}^{\frac{m}{2}-1} C(j) x^{m-2j} y^{m+2j-2}\bigg) + y^{2m-2} \right].
\hskip 2.3cm 
\end{eqnarray}

\noindent Replacing (\ref{coefA}), (\ref{coefB}) and (\ref{coefC}) in (\ref{suma-simpli}) we conclude that
\bigskip

\noindent $P_x \,(Q_y P_{xy} - Q_x P_{yy}) + P_y \,(Q_x P_{xy} - Q_y P_{xx}) =$
{\small
\begin{eqnarray*}
= 2\,k \,(x^2+y^2)^{k-1} m^2 (m-1)
\left[\sum_{j=0}^{\frac{m}{2}-1} \left(
\begin{array}{c}
m-1\\ j\\
\end{array}
\right) x^{2m-2j-2} y^{2j}  \,\,+ \right.\hskip 0.9cm\\
\left. \hskip 2.3cm +
\left(
\begin{array}{c}
m-1\\ \frac{m}{2}\\
\end{array}
\right)
x^{m-2} y^{m} + \sum_{r=\frac{m}{2}+1}^{m-2}
\left(
\begin{array}{c}
m-1\\ r\\
\end{array}
\right)
x^{2m-2r-2} y^{2r}+ \,y^{2m-2}\right].
\end{eqnarray*}}

\noindent Putting together all the terms lying inside of the square brackets we have
\bigskip

\noindent $P_x \,(Q_y P_{xy} - Q_x P_{yy}) + P_y \,(Q_x P_{xy} - Q_y P_{xx}) =$
{\small
\begin{eqnarray*}
\hskip 3.1cm = 2\,k \,(x^2+y^2)^{k-1} m^2 (m-1)\left[
\sum_{j=0}^{m-1}
\left(
\begin{array}{c}
m-1\\ j\\
\end{array}
\right)
x^{2m-2j-2} y^{2j}\right].
\end{eqnarray*}}

\noindent Observe that the expression exposed inside of the square brackets is
the binomial $\,(x^2+y^2)^{m-1}.$ So, finally
$$P_x(Q_y P_{xy} - Q_x P_{yy})+ P_y(Q_x P_{xy} - Q_y P_{xx}) =
2\,k \, m^2 (m-1)
\left(x^2+y^2\right)^{k+m-2}.
\Box$$

\begin{theorem} \label{conjetura}
The Arnold's conjecture is true. In fact, the number of connected components of $Hyp(n)$
is at least $\left[\frac{n-1}{2}\right]$.
\end{theorem}

\noindent {\bf Proof.}
Proposition \ref{polhyper} asserts that the homogeneous polynomial $f^{m+2k} = P^m Q^{2k},\,$ where
$\,k\geq 1,\ m>max \{2,k\}$ is hyperbolic, meanwhile Proposition $\ref{Qelliptic}$ ensures that the polynomial
$Q^{2k}$ is elliptic. Moreover, Theorem \ref{condPQpar} implies that they
satisfy inequality $(\ref{uno})$.
So, by Corollary \ref{coroindices} we conclude that $\, i_0(II_{f^{m+2k}}) = \frac{2-m}{2}$.
Let $\, n\geq 3$ be a natural number.
Now, we shall determine the number of pairs
$\, (k, m)\in \mathbb N\times \mathbb N \,$ such that
$\, k\geq 1, \,m>max \{2,k\}$ and $\,2k+m=n$. (Table \ref{tabla2})
\begin{itemize}
\item If $n$ is even, the set of pairs is $\,\{(k, n-2k) :  k\geq 1,\ m>max \{2,k\} \} = \{(k, n-2k) : k=1, \cdots ,\frac{n}{2}-2\}$.
Moreover, since  $\,i_0\left(II_{f^{m+2k}}\right) = \frac{2-m}{2} = k+1-\frac{n}{2}$,
then each one of these polynomials belongs to different connected component of $Hyp(n)$.
Adding the connected component determined by  the polynomial $P^n$ of degree $n$, we conclude that the number of connected components
of $Hyp(n)$ is at least $\,\frac{n}{2}-1$.
\item If $n$ is odd, the set of pairs is $\,\{(k, n-2k) :  k\geq 1,\ m>max \{2,k\}\} =
\{(k, n-2k) : k=1, \cdots ,\frac{n}{2}-\frac{3}{2}\}$. Moreover, since
$\,i_0\left(II_{f^{m+2k}}\right) = k+1-\frac{n}{2}$,
 each one of these polynomials belongs to different connected component of $Hyp(n)$.
Adding the connected component determined by the polynomial $P^n$ of degree $n$, we conclude that the number of connected components of $Hyp(n)$ is at least $\,\frac{n-1}{2}$. $\hfill\Box$
\end{itemize}

\begin{table}[htb]
\begin{center}
{\small
\begin{tabular}{ | c | c | c | c | c |}\hline
\multirow{2}{*}{$n$ = deg of $f^{m+2k}$}&\multirow{2}{*}{k}&
\multirow{2}{*}{m}&\multirow{2}{*}{$i_0(II_{f^{m+2k}})=\frac{2-m}{2}
$}&Low bound for the \\
 & & & &number of components \\ \hline\hline
3  & 0 & 3 & $-{1}/{2}$ & 1\\ \hline\hline
4  & 0 & 4 & $-{1}$ & 1 \\ \hline\hline
\multirow{2}{*}  {5}  & 0 & 5 & $-{3}/{2}$ & \multirow{2}{*} {2} \\ \cline{2-4}
& 1 & 3 & $-{1}/{2}$ & \\ \hline\hline
\multirow{2}{*}  {6}  & 0 & 6 & $-{2}$ & \multirow{2}{*} {2} \\ \cline{2-4}
& 1 & 4 & $-{1}$ & \\ \hline\hline
\multirow{3}{*}  {7}  & 0 & 7 & $-{5}/{2}$ & \multirow{3}{*} {3} \\ \cline{2-4}
& 1 & 5 & $-{3}/{2}$ & \\ \cline{2-4}
& 2 & 3 & $-{1}/{2}$ & \\ \hline\hline
\multirow{3}{*}  {8}  & 0 & 8 & $-{3}$ & \multirow{3}{*} {3}\\ \cline{2-4}
& 1 & 6 & $-{2}$ & \\ \cline{2-4}
& 2 & 4 & $-{1}$ & \\ \hline
\end{tabular}}
\end{center}
\caption{Hyperbolic homogeneous polynomials up to degree 8.}
 \label{tabla2}
\end{table}

\begin{remark}
The number of connected components for degrees $n=3,4$ and $5,$ was determined with a different approach in  {\rm \cite{geoast}}. 
\end{remark}

Let us provide a qualitative description of the foliation of the field of 
asymptotic lines of the polynomials $f^{2k+m}$, see \cite{hartman}, p. 161.

\begin{corollary}\label{foliations}
The foliation of the field of asymptotic lines of the polyno\-mials
$\,f^{2k+m},\ k \geq 0,\ m \geq 3$ on $\mathbb R^2$
has only one singularity at the origin where $m$ separatrices pass 
through dividing the plane in $m$  
hyperbolic sectors.
\end{corollary}

\noindent {\bf Proof.} Since the second fundamental forms of  
$f^{2k+m}$ and $P^m$ are hyperbolic isotopic their fields of asymptotic 
lines are topologically equivalent. That is, it is enough to describe 
the foliation corresponding to $P^m$.
Considering the natural identification of $\mathbb R^2$ with the complex 
plane, it is easy to define a hyperbolic isotopy between the second fundamental 
form of $P^m$, when $m$ is even, and the hyperbolic quadratic differential form  
Im$(z^{m-2})dz^2$, where $z=x+iy$, $dz=dx + idy$ and Im$(z^{m-2})dz^2$ 
means the imaginary part of the quadratic differential form, 
described by Hopf in \cite{hopf}. If $m$ is odd, the proof follows
by noting that the second
fundamental form of the polynomial $P^m$ composed with the reflection
$T :\mathbb R^2 \rightarrow \mathbb R^2,\, T(u,v)=(v,u)$ is equal to the 
quadratic form $\,m(m-1)$ Im$(z^{m-2})dz^2$.
\hfill$\Box$

\medskip



\noindent {\bf Acknowledgments} We would like to thank Criel Merino 
for his proof of equation $(12)$, which inspired us to use 
similar arguments to prove a couple of useful combinatorial equations. We also thank Francisco Larri\'on and Angel Tamariz for nice conversations on this subject.

\noindent Adriana Ortiz Rodr\'{\i}guez, Instituto de Matem\'aticas,
Universidad Nacional Aut\'onoma de M\'exico,
Ciudad Universitaria, M\'exico D.F 04510, M\'exico.
e-mail: aortiz@math.unam.mx

\noindent Federico S\'anchez-Bringas, Departamento de Matem\'aticas,
Facultad de Ciencias, Universidad Nacional Aut\'onoma de M\'exico,
Ciudad Universitaria, M\'exico D.F 04510, M\'exico.
e-mail: sanchez@servidor.unam.mx

\end{document}